\documentclass[11pt]{article}
\def\bn{\hbox{\it I\hskip -2pt N}}
\def\bz{\hbox{\it Z\hskip -4pt Z}}
\def\bq{\hbox{\it l\hskip -5.5pt Q}}

\def\demo{\noindent{\bf Proof.}}
\newtheorem{theorem}{Theorem}[section]

\newtheorem{remark}[theorem]{Remark}

\newtheorem{definition}[theorem]{Definition}
\newtheorem{proposition}[theorem]{Proposition}

\begin{document}
\begin{center}
\uppercase{{\bf On the geometry of complete intersection\\ toric varieties}}
\end{center}
\advance\baselineskip-3pt
\vspace{2\baselineskip}
\begin{center}
{\sc Nickolas J. Michelacakis and Apostolos Thoma}
\end{center}
\vskip.5truecm\noindent
\vskip.5truecm\noindent
\begin{abstract}
In this paper we give a geometric characterization of the cones of toric varieties that 
are complete intersections. In particular, we prove that the class of complete  
intersection cones is the smallest class of cones which is closed under direct sum 
and contains all simplex cones.  Further, we show that the number of the extreme  
rays of such a cone, which is less than or equal to $2n-2$, is exactly $2n-2$ if and only 
if the cone is a bipyramidal cone, where $n>1$ is the dimension of the cone. Finally, we characterize all toric varieties  
whose associated cones are complete intersection cones. 
\end{abstract}

\section{Introduction} 
  
Motivation for this article has been provided by several results suggesting  
the existence of a strong relation between the geometry of the cone and the 
property of a toric variety to be a complete intersection.   
J. C. Rosales and P. A. Garcia-Sanchez proved in \cite{R-G} that the cone of a 3-dimensional  
affine toric variety contains no more than $4$ extreme rays. In 1997 K.G. Fischer, W. Morris   
and J. Shapiro \cite{F-M-S} generalized this result by proving that the cone 
of an n-dimensional affine toric variety contains no more than $2n-2$ extreme rays. 
This result was also proved, independently, by M. Katzman \cite{K} in a graph theory context. 
  
The purpose of this article is to investigate this relation. We provide a contrapositive  
geometric criterion enabling one to decide when an affine toric variety is not a complete 
intersection. 
  
Let $A=\{{\bf a}_i \,|\, 1\leq i \leq m\}$ be a subset of $\bz^n$ such that the  
 semigroup generated by $A$, $\bn A=  
\{\sum _{i=1}^mn_i{\bf a}_i \,|\, {\bf a}_i\in A,\ n_i\in \bn \}$,  
is an  affine semigroup. An {\em affine semigroup} $S$ is a finitely generated 
subsemigroup of $\bz^n$ with no invertible elements, that means $S\cap (-S)=\{{\bf 0}\}$. 
To the set $A$ we associate the 
 {\em toric ideal} $I_A$ which is the kernel of the $K$-algebra homomorphism 
$$\phi :K[x_1,\ldots ,x_m] \rightarrow K[t_1,\ldots ,t_n,t_1^{-1},\ldots ,t_n^{-1}]$$  
given by 
$$\phi (x_i)={\bf t}^{{\bf a}_i}=t_1^{a_{i,1}}\dots t_n^{a_{i,n}} \ \mbox{for all}\ i=1,\ldots ,m\, ,$$ 
where ${\bf a}_i=(a_{i,1},\ldots ,a_{i,n})$. 
The set of zeroes  of $I_A$, $V(I_A)\subset K^m$  is an {\em affine toric variety} in the sense of \cite{St}, since  
we do not require normality. To the toric variety $V(I_A)$ we associate a strongly convex rational  
polyhedral cone 
\[  
	\sigma =pos_{\bq }(A):=\{l_1{\bf a}_1+ \cdots +l_m{\bf a}_m\, |\, l_i\in \bq^+\}:=\bq^+ A,  
\]  
where $\bq^+$ denotes the set of nonnegative rationals.  
  
An affine toric variety $V(I_A)$ (resp. an affine semigroup $\bn A$)  
is called {\em complete intersection } 
if and only if the ideal $I_A$ is a 
complete intersection, i.e. the minimal number of generators of $I_A$ equals  
the height of $I_A$.  
For $A=\{{\bf a}_1,\cdots ,{\bf a}_m\}\subset \bz ^n$ and $E\subset \{ 1,\dots ,m\}$, we denote  
by $A^E$ the set $\{{\bf a}_i| i\in E\}$. 
\begin{definition} 
	Let $E_1$, $E_2$ be two nonempty subsets of $\{1, \ldots ,m\}$ such that  
	$E_1\cup E_2=\{1, \ldots ,m\}$ and $E_1\cap E_2=\emptyset $. 
	The semigroup $\bn A$ is called the gluing of the semigroups 
	$\bn A^{E_1}$ and $\bn A^{E_2}$ if there is a nonzero 
	${\bf a}\in \bn A^{E_1}\cap \bn A^{E_2}$ such that 
	$\bz{\bf a}=\bz A^{E_1}\cap \bz A^{E_2}$.  
\end{definition} 
 
Complete intersection affine semigroups or equivalently complete intersection  
toric varieties have been completely characterized in \cite{F-M-S} via the use of semigroup gluings. 
 
\begin{theorem}[Fischer-Morris-Shapiro]   
	Let $\bn A$ be an affine semigroup that is not a free abelian semigroup.  
	Then $\bn A$ is a complete intersection if and only if $\bn A$ is  
	the gluing of the semigroups 
	$\bn A^{E_1}$ and $\bn A^{E_2}$ and both $\bn A^{E_1}$, $\bn A^{E_2}$  
	are complete intersections.  
\end{theorem} 
As an immediate corollary, one gets:  

\begin{theorem} The class of complete intersection  affine semigroups is the smallest class of  
	affine semigroups closed under the operation of gluing and containing all free affine semigroups.  
\end{theorem}  
In this article, we use this characterization of complete intersection affine semigroups 
to provide a geometric characterization of complete intersection cones. A cone $\sigma $ is called a 
{\em complete intersection cone} if there exists a complete intersection toric  
variety $V(I_A)$  such that its associated cone $\sigma_A$ is equal to $\sigma $.  
More specifically, in section 2, we present the notion of the direct sum of cones and prove its basic properties.  
In section 3, we characterize complete intersection cones by proving that the class of complete intersection cones  
is the smallest class of cones which is closed under direct sum and contains all simplex cones. 
Further, we show that the number of the extreme  
rays of such a cone, which is less than or equal to $2n-2$, is exactly $2n-2$ if and only 
if the cone is a bipyramidal cone. 
In section 4, we characterize all toric varieties  
whose corresponding cones are complete intersection cones. 
  
\section{Direct sum of cones} 
  
Let $B\subset\bq^n$.  
The rational linear hull of the set $B$ is the set  
\[  
	\mbox{span}_{\bq}(B):=\{\sum _{i=1}^lq_i{\bf b}_i\, |\, q_i\in \bq ,\, {\bf b}_i\in B\}:=\bq B,  
\]  
while the rational polyhedral cone of $B$  is the set 
\[  
	pos_{\bq }(B):=\{\sum _{i=1}^lq_i{\bf b}_i\\, |\, q_i\in \bq^+,\, {\bf b}_i\in B\}:=\bq^+ B.  
\]  
The relative interior  
$\mbox{relint}(\sigma)$ of a cone $\sigma $ is the usual topological interior of $\sigma$ considered  
as a subset of $\mbox{span}_{\bq}(\sigma)$.  
A cone $\sigma $ is {\em strongly convex} if $\sigma \cap (-\sigma )=\{ {\bf 0} \}$.  
A face $F$ of ${\sigma}$ is any set of the form 
	$$F= \sigma \cap \{{\bf x}\in \bq ^n: {\bf c}\cdot {\bf x} = 0\}, $$  
where ${\bf c}\in \bq ^n$ and  ${\bf c}\cdot {\bf x} \ge 0$ 
for all points ${\bf x}\in \sigma$. Note that   
${\bf c}\cdot {\bf x}$ denotes the dot product in $\bq ^n$. We say that the 
vector ${\bf c}$ defines the face $F$. 
Faces of dimension one are called {\em extreme rays}. An extreme ray will be denoted by 
$\bq ^+{\bf r}$ for a nonzero vector, ${\bf r}$, belonging to it. The set of 
extreme rays of a cone $\sigma $ 
will be denoted by $\mbox{Rays}(\sigma )$.  If the number of the extreme 
rays $k$ coincides with the dimension of the cone (i.e. the extreme rays are linearly  
independent), the cone is called a $k$-{\em simplex cone}.  
  
\begin{definition} 
	Let $\sigma _1$, $\sigma _2$ strongly convex rational polyhedral cones in $\bq ^n$ 
	with $\mbox{span}_{\bq}(\sigma _1)\cap \mbox{span}_{\bq}(\sigma _2)=\bq {\bf a}$ 
	and ${\bf a} \in \sigma _1\cap \sigma _2$. We call the 
	rational polyhedral cone $\sigma _1  \oplus   \sigma _2:=pos_{\bq }(\sigma _1\cup \sigma _2)$  
	the direct sum of $\sigma _1,   \sigma _2$ along ${\bf a}$.  
	If ${\bf a}$ belongs to an extreme ray $\bq^+{\bf r}$ of either $\sigma _1$, $\sigma _2$  or both  
	we call the direct sum, direct sum of external type and we 
	denote it $\sigma _1  \oplus_{\bf r}  \sigma _2$. 
	Otherwise we call the direct sum, direct sum of internal 
	type and we denote it $\sigma _1  \oplus_{\bf a}  \sigma _2$. 
\end{definition}  
  
We shall make use of the notation $\sigma _1\oplus\sigma _2$ to indicate  
the direct sum of two cones $\sigma_1$ and $\sigma_2$ in general, regardless of type.  
Observe that given two cones $\sigma _1=pos_{\bq }({\bf r}_{11}, \dots ,{\bf r}_{1t})$ and  
$\sigma _2=pos_{\bq }({\bf r}_{21}, \dots ,{\bf r}_{2s})$ then, $\sigma _1  \oplus  \sigma _2= 
pos_{\bq }({\bf r}_{11}, \dots ,{\bf r}_{1t},{\bf r}_{21}, \dots ,{\bf r}_{2s} )$. Thus, 
every element of $\sigma $ can be written as a sum of an element in $\sigma _1$ and an element 
in $\sigma _2$, albeit not necessarily in a unique way.\\

\begin{proposition} \label{1}  
	The rational polyhedral cones $\sigma _1$, $\sigma _2$ are 
	strongly convex if and only if  
	$\sigma _1  \oplus  \sigma _2$ is strongly convex. 
\end{proposition}  
  
\demo Assume that both $\sigma _1$, $\sigma _2$ are strongly convex and set  
$\sigma =\sigma _1  \oplus  \sigma _2$. Let ${\bf x}\in \sigma \cap (-\sigma )$. Then,  
${\bf x}={\bf x}_1+{\bf x}_2=-{\bf y}_1-{\bf y}_2$, where ${\bf x}_1, {\bf y}_1\in \sigma_1$ 
and ${\bf x}_2, {\bf y}_2\in \sigma_2$. Set ${\bf w}:={\bf x}_1 + {\bf y}_1=-{\bf x}_2-{\bf y}_2$, 
then ${\bf w}\in \mbox{span}_{\bq}(\sigma _1)\cap \mbox{span}_{\bq}(\sigma _2)=\bq {\bf a}$,  
i.e. ${\bf w}=q{\bf a}$.
From ${\bf x}_1 + {\bf y}_1=q{\bf a}$ we conclude that $q\geq 0$ 
and from $-{\bf x}_2-{\bf y}_2=q{\bf a}$ 
we conclude that $q\leq 0$. Thus, $q=0$ which implies that ${\bf x}_1 =- {\bf y}_1$ and  
${\bf x}_2 =- {\bf y}_2$. Because both $\sigma _1$, $\sigma _2$ are strongly convex  
${\bf x}_1={\bf x}_2={\bf 0}$, therefore, ${\bf x}={\bf 0}$ and 
$\sigma _1  \oplus   \sigma _2$ is strongly convex. \\ 
The reverse inclusion follows immediately since $\sigma _1,\sigma _2 \subset \sigma _1  
\oplus  \sigma _2$.  
  
\begin{proposition} \label{2} 
	Simplex cones are strongly convex. 
\end{proposition} 
  
\demo Let $\sigma =pos_{\bq }({\bf r_1},\cdots ,{\bf r_s})$ such that  
${\bf r_1},\cdots ,{\bf r_s}$ are linearly independent. Let ${\bf x}\in \sigma \cap (-\sigma )$. 
Then ${\bf x}=a_1{\bf r}_1+\cdots +a_s{\bf r}_s=-b_1{\bf r}_1+\cdots -b_s{\bf r}_s$, where 
$a_i, b_i$ are nonnegative rationals. Then, $(a_1+b_1){\bf r}_1+\cdots +(a_s+b_s){\bf r}_s={\bf 0}$ 
implies that $a_i+b_i=0$ for all $i$. Therefore $a_i=0=b_i$, which 
means that ${\bf x}={\bf 0}$.  
  
When the cones are strongly convex we can find a hyperplane intersecting the cones $\sigma $,   
$\sigma _1$ and $\sigma _2$ at  convex polytopes $P$, $P_1$ and $P_2$ 
\cite{Ewald}. 
Translating the definition of the direct sum of cones to the polytopes, we have 
the notion of direct sum of polytopes. 
In case ${\bf a}$ is in the relative 
interiors of both $\sigma _1$ and $\sigma _2$ 
then the polytope $P$ is the usual  direct sum (or free sum) of $P_1$ and $P_2$ \cite{Ewald,R-Z}. 
It is, however, important to point out that 
in our case ${\bf a}$ may belong to the 
relative interior of any k-dimensional face, $k\geq 2$, in the direct sum of internal 
type, or to an extreme ray, in the direct sum of external type.  
  
\begin{theorem} \label{3} 
	Let $\sigma =\sigma _1  \oplus   \sigma _2$ with 
	$dim(\sigma _1 )=n_1$ and $dim(\sigma _2)=n_2$. 
	Let $k_1$ be the number of extreme rays of $\sigma _1$ and $k_2$ be the 
	number of extreme rays of $\sigma _2$, 
	then $\sigma _1  \oplus   \sigma _2$ has dimension $n_1+n_2-1$, $k_1+k_2$ 
	extreme rays if the direct 
	sum is of internal type 
	and $k_1+k_2-1$ if the direct 
	sum is of external type. 
\end{theorem} 
  
\demo The formula for the dimension follows immediately from the hypothesis 
since $\mbox{span}_{\bq}(\sigma _1)\cap \mbox{span}_{\bq}(\sigma _2)=\bq {\bf a}$.  
We show first that 
\begin{equation}\label{eq1}  
	\mbox{Rays}(\sigma )\subset \mbox{Rays}(\sigma _1)\cup \mbox{Rays}(\sigma _2). 
\end{equation}  
To this end,  
let $\bq^+{\bf r}$ be an extreme ray of $\sigma$ then, ${\bf r}={\bf x}_1+{\bf x}_2$, with 
${\bf x}_1\in \sigma _1$, ${\bf x}_2\in \sigma _2$. Let further, ${\bf c}\in \mbox{span}_{\bq}(\sigma ):=  
\bq ^{n_1+n_2-1}$ 
be a vector defining the extreme ray $\bq^+{\bf r}$. Then, ${\bf c}\cdot {\bf x}=0$ for all 
${\bf x}\in\bq^+{\bf r}$, while ${\bf c}\cdot {\bf x}>0$ for all other ${\bf x}\in\sigma$. 
This implies  
that ${\bf c}\cdot {\bf r}=0$ while ${\bf c}\cdot {\bf x}_1\geq 0$ and  
${\bf c}\cdot {\bf x}_2\geq 0$, therefore they should both belong to $\bq _+{\bf r}$.  
Since ${\bf r}\not ={\bf 0}$ at least one of ${\bf x}_1, {\bf x}_2$ is nonzero. 
Then every nonzero vector of ${\bf x}_1, {\bf x}_2$ defines an extreme ray of $\sigma _1$,   
$ \sigma _2$ respectively, i.e. 
every ray of $\sigma $ is either an extreme ray of $\sigma _1$ or of $ \sigma _2$ 
proving (\ref{eq1}).\\ 
We now look at the extreme rays of $\sigma_1, \sigma_2$. \\  
Choose a basis of $\mbox{span}_{\bq }(\sigma _1)$ such that the 
$n_1$-th element of the basis is either ${\bf a}$, if $\sigma$ is of internal type,  
or ${\bf r}$, if $\sigma$ is of external type, and extend it to a basis of  
$\mbox{span}_{\bq }(\sigma )=\bq ^{n_1+n_2-1}$.  
Then, for ${\bf x}\in \sigma _1$, we have that $({\bf x})_j=0$  
for all $n_1<j\leq n_1+n_2-1$ while for ${\bf y}\in \sigma _2$,  
we have that $({\bf y})_j=0$ for all $1\leq j< n_1$.\\  
Let $\bq ^+{\bf r_1}$ be any ray of $\sigma _1$, if $\sigma$ is 
of internal type, and any ray of $\sigma _1$ different from $\bq ^+{\bf r}$, 
if $\sigma$ is of external type. 
Since $\bq ^+{\bf r_1}$ is an extreme ray of $\sigma _1$ 
there exist ${\bf c}_1\in \bq ^{n_1+n_2-1}$  
such that ${\bf c}_1\cdot {\bf r}_1=0$ and ${\bf c}_1\cdot {\bf x}>0$  
for all ${\bf x}\not \in \bq ^+{\bf r_1}$.  
The cone $\sigma _2$ is strongly convex thus, $\{{\bf 0}\}$ is a face 
of $\sigma _2$. So, there exists a ${\bf c}_2\in \bq ^{n_1+n_2-1}$  
such that ${\bf c}_2\cdot {\bf y}>0$ for all  
${\bf 0}\not= {\bf y}\in \sigma _2$. We may assume that $({\bf c}_1)_j=0$  
for all $n_1<j\leq n_1+n_2-1$ and $({\bf c}_2)_j=0$ for all $1\leq j< n_1$. 
Note that ${\bf a}\not \in \bq ^+{\bf r}_1$, if $\sigma$ of internal type, and 
${\bf r}\not \in \bq ^+{\bf r}_1$, if $\sigma$ external type,  
therefore, 
${\bf c}_1\cdot {\bf a}>0$ and ${\bf c}_2\cdot {\bf a}>0$ (resp. ${\bf c}_1\cdot {\bf r}>0$ and 
${\bf c}_2\cdot {\bf r}>0$)  
 which means that $({\bf c}_1)_{n_1}>0$ and $({\bf c}_2)_{n_1}>0$.  
 Then the vector   	  
$$	({\bf c})_i= 		\left\{ 			  
	\begin{array}{lll} 			 	  
	({\bf c}_2)_{n_1}({\bf c}_1)_i & i=1,\ldots ,n_1-1 \\ 			 	  
	({\bf c}_2)_{n_1}({\bf c}_1)_{n_1} & i=n_1	\\ 		  
	({\bf c}_1)_{n_1}({\bf c}_2)_i & i=n_1+1,\ldots ,n_1+n_2-1 			  
	\end{array} 		  
	\right.	$$  
defines ${\bf r}_1$ as an extreme ray of $\sigma $. A similar argument shows that 
the corresponding statement concerning the extreme rays of $\sigma _2$ holds true.  
Therefore, every extreme ray of $Rays(\sigma _1)\cup Rays(\sigma _2)$,  
different from $\bq ^+{\bf r}$,  
is an extreme ray of $\sigma $. Thus, if  
$\sigma =\sigma _1  \oplus _a   \sigma _2$ is of internal type,  
we have proved that $\mbox{Rays}(\sigma )= 
\mbox{Rays}(\sigma _1)\cup \mbox{Rays}(\sigma _2)$ and since  
$\mbox{Rays}(\sigma _1), \mbox{Rays}(\sigma _2)$ have no common elements  
we conclude that $\sigma =\sigma _1  \oplus _a   \sigma _2 $ has $k_1+k_2$ extreme rays. \\  
If $\sigma$ is of external type, we need to check what happens with 
 $\bq ^+{\bf r}$. There are two cases:\\  
 i) $\bq ^+{\bf r}$ is an extreme ray of both 
$\sigma _1, \sigma _2$. Then there exist ${\bf c}_1, {\bf c}_2\in \bq ^{n_1+n_2-1}$ such that 
${\bf c}_1\cdot {\bf r}=0$ and ${\bf c}_1\cdot {\bf x}>0$ for all ${\bf x}\in \sigma _1$ and  
not in $\bq ^+{\bf r}$,  
and ${\bf c}_2\cdot {\bf r}=0$ and ${\bf c}_2\cdot {\bf y}>0$ for  
all ${\bf y}\in \sigma _2$ and not in $\bq ^+{\bf r}$.
From ${\bf c}_1\cdot {\bf r}=0$, ${\bf c}_2\cdot {\bf r}=0$ we get $({\bf c}_1)_{n_1}=({\bf c}_2)_{n_1}=0$  
since the vector ${\bf r}$ is the $n_1$-th element in the basis of $\bq ^{n_1+n_2-1}$.  
Then the vector   	  
$$	({\bf c})_i= 		\left\{ 			  
	\begin{array}{lll} 			 	  
	({\bf c}_1)_i & i=1,\ldots ,n_1-1 \\ 			 	  
	({\bf c}_1)_{n_1}=({\bf c}_2)_{n_1}=0 & i=n_1	\\ 		  
	({\bf c}_2)_i & i=n_1+1,\ldots ,n_1+n_2-1 			  
	\end{array} 		  
	\right.	$$  
defines $\bq ^+{\bf r}$ as an extreme ray of $\sigma $. The extreme rays in this case are all  
rays in $\mbox{Rays}(\sigma _1)\cup \mbox{Rays}(\sigma _2)$, but since $\bq ^+{\bf r}$  
is counted twice   
we have a total of $k_1+k_2-1$ rays.\\  
ii) $\bq ^+ {\bf r}$ is an extreme ray of one of them; we may, without loss of generality,  
assume of $\sigma _1$.  
We claim that $\bq ^+ {\bf r}$ cannot be an extreme ray of 
$\sigma $. If it were, there would exist ${\bf c}\in \bq ^{n_1+n_2-1}$ 
such that ${\bf c}\cdot {\bf r}=0$ and ${\bf c}\cdot {\bf x}>0$ 
for all ${\bf x}\in \sigma $ not in $\bq ^+ {\bf r}$.  
But then, $\bq ^+ {\bf r}\subset \sigma _2$, ${\bf c}\cdot {\bf r}=0$ and 
${\bf c}\cdot {\bf x}>0$ for every 
${\bf x}\in \sigma _2$ not in $\bq ^+ {\bf r}$ forcing 
$\bq ^+ {\bf r}$ to be an extreme ray of $\sigma _2$, a 
contradiction. 
Thus, in this case the extreme rays of $\sigma =\sigma _1  \oplus _r   \sigma _2$ are all the extreme rays  
of $\mbox{Rays}(\sigma _1)\cup \mbox{Rays}(\sigma _2)$ except $\bq ^+ {\bf r}$ which makes 
again a total of $k_1+k_2-1$ rays.

\begin{remark} \label{4}  
	If $dim(\sigma _2)=1$ then, $\sigma_1\oplus\sigma_2$ is always of external type and $\sigma_1\oplus_{\bf r}\sigma_2=\sigma _1$. 
\end{remark}

\section{Complete intersection cones} 
  
We fix $A=\{{\bf a}_1,\cdots ,{\bf a}_m\}\subset \bz ^n$ and $E\subset \{ 1,\dots ,m\}$  
and denote by  $\sigma _E$ the cone $pos_{\bq}({\bf a}_i| i\in E)=\bq A^{E}$.  
We have the following proposition:

\begin{proposition} \label{5}  
	If the  semigroup $\bn A$ is  
	the gluing of the semigroups $\bn A^{E_1}$ and $\bn A^{E_2}$  then  
	$\sigma =\sigma _{E_1}  \oplus  \sigma _{E_2}$.  
\end{proposition}  
  
\demo Suppose that $\bn A$ is the gluing of the 
semigroups $\bn A^{E_1}$ and $\bn A^{E_2}$. 
Then there exist an ${\bf a}$ such that $\bz{\bf a}=\bz A^{E_1}\cap \bz A^{E_2}$, therefore  
$\mbox{span}_{\bq}(\sigma _{E_1})\cap \mbox{span}_{\bq}(\sigma _{E_2})=
\bq  A^{E_1}\cap \bq A^{E_2}=\bq {\bf a}$. 
From ${\bf a}\in \bn A^{E_1}\cap \bn A^{E_2}$ we get that ${\bf a} 
\in \bq ^+ A^{E_1}\cap \bq ^+ A^{E_2}= 
\sigma _{E_1}\cap \sigma _{E_2}$.

\begin{definition} 
	A cone $\sigma $ is called a complete intersection cone if there 
	exist a complete intersection 
	toric ideal whose associated cone is $\sigma $. 
\end{definition}  
  
\begin{theorem} 
	The class of complete intersection cones is the smallest class of 
	cones closed under direct sum containing all rational simplex cones. 
\end{theorem} 
  
\demo By propositions \ref{1} and  \ref{2} all cones 
in the above class  are strongly convex. 
Let $\sigma $ be a rational simplex cone. Then by choosing an integer vector 
in each extreme ray we have a set $A\subset\bz ^n$ of linearly independent 
vectors over $\bq$. Equivalently, 
the ideal $I_A$ is the zero ideal and thus a complete intersection, and so is $\bn A$.  
Suppose that $\sigma =\sigma_1\oplus\sigma_2$ is a direct sum of 
two complete intersection cones $\sigma _1$ and $\sigma _2$.  
This means that there exist two sets of integer vectors such that 
$\sigma _1  =pos_{\bq }(A_1)$, $\sigma _2 =pos_{\bq }(A_2)$ and $\bn A_1$, 
$\bn A_2$ are complete intersections.  
Since $\sigma =\sigma _1  \oplus    \sigma _2$ we have that 
$\mbox{span}_{\bq} (A_1)\cap \mbox{span}_{\bq} (A_2)$ is a $\bq $-vector space of 
dimension 1 therefore  
there exist an integer vector ${\bf a}\in \sigma$ such that $\bz A_1\cap \bz A_2  
=\bz {\bf a}$.  Since  
$\bz A_i\cap \bz {\bf a}=\bz {\bf a},\, i=1,2$, the semigroups $\bn A_i\cap \bn {\bf a}
=<n_{i,1}{\bf a},\cdots ,n_{i,s_i}{\bf a}>$  
with $g.c.d.(n_{i,1},\cdots ,n_{i,s_i})=1,\, i=1,2$, i.e. they are isomorphic to numerical semigroups. 
Therefore, for all sufficiently large $n$, $n{\bf a}\in \bn A_i,\, i=1,2$.\\ 
Let ${\bf a}=(a_1,\ldots ,a_n)$,  
 $g=g.c.d.(a_1,\ldots ,a_n)$ and $\tau $ be any positive integer prime with $g$. 
We  claim that if $\bz A_1\cap \bz A_2  
=\bz {\bf a}$ then, $\bz \tau A_1\cap \bz A_2=\bz \tau {\bf a}$. To prove it, 
if ${\bf b}\in \bz \tau A_1\cap \bz A_2\subset \bz A_1\cap \bz A_2  
=\bz {\bf a}$ then, ${\bf b}=t{\bf a}$ for an integer $t$. 
Since $t{\bf a}\in \bz \tau A_1$,  
each coordinate of $t{\bf a}$ is a multiple of $\tau $.  
But $(\tau,g)=1$ thus, $t$ is a multiple of $\tau $, proving the claim.\\ 
Choose $\mu $ and $\tau $ positive integers such that they are prime to each other 
and to $g$ and such that  
$\mu {\bf a}\in \bn A_1$ and $\tau {\bf a}\in \bn A_2$. Then,  
by the previous argument applied twice, we  
have that $\bz \tau A_1\cap \bz \mu A_2=\bz \tau \mu {\bf a}$. But also, 
$\tau \mu {\bf a}\in \bn \tau A_1\cap \bn \mu A_2$ therefore, the toric ideal 
for $A=\tau A_1 \cup \mu A_2$ is  
a complete intersection since it is the gluing of $\bn \tau A_1$, 
$\bn \mu A_2$ and both of them are  
complete intersections.\\  
The other direction follows by translating Theorem 1.3 from semigroup gluings to direct sums  
of cones using Proposition 3.1.

\begin{definition}  
	We call general bipyramidal cones the elements of the smallest class 
	of cones which is closed under the operation of direct sum of 
	internal type and contains all  2-simplex cones. 
\end{definition}  
  
In \cite{F-M-S}, K. G. Fisher, W. Morris and J. Shapiro proved that the number of extreme  
rays of a complete intersection cone of dimension $n>1$ is less 
than or equal to $2n-2$. In order to show that the bound $2n-2$, 
for the number of extreme rays, is best possible they presented a family of examples,  
one for each dimension. All cones in these examples were bipyramidal.  
We provide next, a direct proof of this result  
which enables us to show that the class of complete intersection cones featuring 
the maximum number, $2n-2$, of extreme rays is precisely the class of bipyramidal cones.  
  
\begin{theorem} The number of extreme rays of a complete intersection cone of dimension $n>1$ is less 
than or equal to $2n-2$. It is 
exactly $2n-2$ if and only if it is a general bipyramidal cone. 
\end{theorem}  
 
\demo Let $\sigma =pos_{\bq }(A)$ be a complete intersection 
cone with 
 $A=\{{\bf a}_1,\dots ,{\bf a}_m\}$ such that $\bn A$ is a complete 
intersection. 
Suppose that $dim({\sigma })=n$ 
and $\sigma $ has $k$ extreme rays. The height of $I_A$ is $m-n$ therefore 
the ideal $I_A$ has a minimal set of 
 generators consisting of $m-n$ binomials.\\ 
A partition ${\cal{J}}$ of $A$ is a set of pairwise disjoint nonempty 
subsets covering $A$. We say that partition ${\cal{J}}_2$ refines partition ${\cal{J}}_1$ 
and write ${\cal{J}}_2 < 
{\cal{J}}_1$, if every set in ${\cal{J}}_1$ is the union of some sets in 
${\cal{J}}_2$. 
According to K. G. Fisher, W. Morris and J. Shapiro \cite{F-M-S} , see also 
Theorems 1.2, 1.3, we have a chain of 
partitions of $A$: 
$${\cal{J}}_{m-n+1} < {\cal{J}}_{m-n} 
< \cdots {\cal{J}}_{2} < {\cal{J}}_{1}$$ 
such that i) ${\cal{J}}_{1}=\{A\}$, \\ ii) for $i=1,\ldots ,m-n$, 
partition ${\cal{J}}_{i}$ is 
obtained from ${\cal{J}}_{i+1}$ by replacing only two sets $E_1, E_2$ of 
the partition 
 ${\cal{J}}_{i+1}$ by their union $E$, while all the rest remain the same; 
the semigroup $\bn E$ is the 
gluing of $\bn E_1$ and $\bn E_2$ 
and\\ 
iii) the vectors in each set of the last partition ${\cal{J}}_{m-n+1}$ are 
linearly independent.\\ 
Note that the partition ${\cal{J}}_{i}$ consists of $i$ sets and 
for any set $S$ in the last partition ${\cal{J}}_{m-n+1}$ we have 
$pos _{\bq }(S)$ is 
a simplex cone. \\ 
From ii) and Proposition \ref{5} we have $pos _{\bq} (E)= 
pos _{\bq} (E_1) \oplus pos _{\bq} (E_2)$, where $pos _{\bq} (E_1) \oplus pos 
_{\bq} (E_2)$ is either of internal type or of external type. The cone $\sigma $ 
is, thus, obtained by taking $m-n$ direct 
sums starting with the $m-n+1$ simplex cones $pos _{\bq} (S)$, where $S\in 
{\cal{J}}_{m-n+1}$. 
Suppose that $a$ of the $m-n$ direct sums are of internal type and $b$ of external type. 
It follows that 
$$a+b=m-n.$$ 
 
If $\sigma$ is a strongly rational polyhedral cone, let
${\cal{R}}(\sigma )$  denote the number of extreme rays 
of $\sigma $. 
For every partition ${\cal{J}}_{i}$ we have a set of cones $\{ 
pos_{\bq}(S) | S\in {\cal{J}}_{i}\}$. 
To any partition we assign the quantity ${\cal{D}}({\cal{J}}_{i})$ which 
is equal to 
$$\sum _{S\in 
{\cal{J}}_{i}}({\cal{R}}(pos_{\bq}(S))-dim_{\bq}(pos_{\bq}(S))).$$ 
The set of cones corresponding to partitions ${\cal{J}}_{i}, {\cal{J}}_{i+1}$,  
differ only in that two cones in ${\cal{J}}_{i+1}$ are replaced by their 
direct sum in ${\cal{J}}_{i}$. Theorem \ref{3} says that 
${\cal{D}}({\cal{J}}_{i})={\cal{D}}({\cal{J}}_{i+1})$ if this is
of external type, while ${\cal{D}}({\cal{J}}_{i})={\cal{D}}({\cal{J}}_{i+1})+1$
if it is of internal type. Since ${\cal{J}}_{1}$ contains just $A$,
or equivalently only the cone $\sigma$, ${\cal{D}}({\cal{J}}_{1})=k-n$, while 
${\cal{D}}({\cal{J}}_{m-n+1})=0$, because every 
cone in ${\cal{J}}_{m-n+1}$ is simplex cone. We conclude that 
$$k-n=a.$$ 
Since in the last partition all the $m$ vectors in $A$ appear as 
extreme rays of the $m+n-1$ 
cones of $pos_{\bq}(S)$, we have 
$\sum _{j=1}^{m+n-1}({\cal{R}}(pos_{\bq}(S_j)))=m$, i.e. 
$$\sum _{j=1}^{m+n-1}({\cal{R}}(pos_{\bq}(S_j))-1)=n-1.$$ 
We may assume that the first $t$ terms in the left part 
of the above equation are nonzero. For the zero terms,
we have that  $({\cal{R}}(pos_{\bq}(S_j)))=1$,  $t<j\leq m+n-1$. 
Remark 2.5 guarantees that the sets of cones of dimension more than 1 in partitions 
${\cal{J}}_{i+1}$, ${\cal{J}}_{i}$ are 
the same and the type of the direct sum involved is external. Thus,
internal type of direct sum may appear 
only amongst the first $t$ cones, $pos_{\bq}(S_j)$. Thus, $$a\leq t-1.$$ 
Further, $\sum _{i=1}^{t}({\cal{R}}(pos_{\bq}(S_i))-1)=n-1$ 
is a nonzero partition of $n-1$ therefore, $$t\leq n-1.$$ 
We conclude that $$a\leq t-1\leq n-2.$$ 
But $k-n=a$ therefore, $k\leq 2n-2$. \\ 
Now, for the second part of the statement of the theorem,
if we have $k=2n-2$ then, $t=n-1$ which implies that each $k_i$ is two or 
one. Actually, we have $n-1$ two's and 
the rest are one's. By remark 2.5 once more, we have that $\sigma $ is 
obtained by taking $n-2$ direct sums of internal type starting with 
$(n-1)$ 2-simplex cones, i.e. $\sigma$ is a general bipyramidal cone.\\ 
For the other direction, let $\sigma $ be a general bipyramidal cone of 
dimension $n$. We claim that 
 the number of extreme rays of $\sigma $ is $2n-2$.  
If $\sigma $ is a 2-simplex cone then, it has $2=2\cdot 2-2$ extreme rays. 
We apply induction assuming the claim true for all bipyramidal cones of dimension 
smaller than $n$, $n\ge 3$. 
Let $\sigma $ be a general bipyramidal cone of dimension $n$, $n\ge 3$, then 
$\sigma$ is the direct sum of internal type of two 
general bipyramidal cones of dimensions $n_1$, $n_2$ smaller than 
$n$. Applying the induction hypothesis on these two cones we see 
that they have $2n_1-2$, $2n_2-2$ extreme 
rays respectively. Now, Theorem 2.4 implies 
that $\sigma $ has $(2n_1-2)+(2n_2-2)=2(n_1+n_2-1)-2=2n-2$ extreme rays, and 
the claim is proved.

\section{Toric varieties with complete intersection cones} 
 
The associated cone a toric variety may be a complete intersection cone 
without the variety  being a complete intersection. Using the results  
 of the previous sections, it is not difficult to determine  
the  vector-sets $A$ so that the cone associated with the toric variety $V(I_A)$ 
(resp. of the affine semigroup $\bn A$) is a  complete intersection cone. In this
description, it turns out that the key concept needed is that of the s-gluing 
of semigroups, a generalization of the 
gluing and p-gluing of semigroups, see relatively \cite{B-M-T,M-T}. 
 
\begin{definition}  
	Let $E_1$, $E_2$ be two nonempty subsets of $\{1, \ldots ,m\}$ such that   
	$E_1\cup E_2=\{1, \ldots ,m\}$ and $E_1\cap E_2=\emptyset $.  
	The semigroup $\bn A$ is called the s-gluing of the semigroups  
	$\bn A^{E_1}$ and $\bn A^{E_2}$ if there exist an integer vector ${\bf a}$ 
	such that $\bz{\bf a}=\bz A^{E_1}\cap \bz A^{E_2}$ 
	and $t{\bf a}\in \bn A^{E_1}\cap \bn A^{E_2}$ 
	for some positive integer $t$. 
\end{definition}  
  
The next proposition  is a generalization of Proposition \ref{5} and makes the relation 
between the operation of s-gluing of semigroups and the operation of taking direct sums 
of cones exact. 
\begin{proposition}  
	The  semigroup $\bn A$ is   
	the s-gluing of the semigroups $\bn A^{E_1}$ and $\bn A^{E_2}$  if and only if  
	$\sigma =\sigma _{E_1}  \oplus  \sigma _{E_2}$.   
\end{proposition} 
 
\demo Suppose that $\bn A$ is the s-gluing of the  
semigroups $\bn A^{E_1}$ and $\bn A^{E_2}$.  
Then there exist an ${\bf a}$ such that $\bz{\bf a}=\bz A^{E_1}\cap \bz A^{E_2}$, therefore   
$\mbox{span}_{\bq}(\sigma _{E_1})\cap \mbox{span}_{\bq}(\sigma_{E_2})=
\bq  A^{E_1}\cap \bq A^{E_2}=\bq {\bf a}$. 
From $t{\bf a}\in \bn A^{E_1}\cap \bn A^{E_2}$ we get that ${\bf a}  
\in \bq ^+ A^{E_1}\cap \bq ^+ A^{E_2}=  
\sigma _{E_1}\cap \sigma _{E_2}$.\\  
For the converse, suppose that $\mbox{span}_{\bq}(\sigma _{E_1})\cap \mbox{span}_{\bq}(\sigma _{E_2})=\bq {\bf a}$ with   
${\bf a} \in \sigma _{E_1}  \cap  \sigma _{E_2}$. Then  
$\bz A^{E_1}\cap \bz A^{E_2}$ is of rank one, and let ${\bf b}$  
be a generator of it such that ${\bf b}=\mu {\bf a}$. Then ${\bf b}=\mu {\bf a}\in   
\sigma _{E_1}  \cap  \sigma _{E_2}=\bq ^+ A^{E_1}\cap \bq ^+ A^{E_2}$.  
Clearing denominators we get that  
$t{\bf b}\in \bn A^{E_1}\cap \bn A^{E_2}$, for some positive integer $t$.

Proposition 4.2 in conjunction with Theorem 3.3 give the following:
 
\begin{theorem} 
	The class of  affine semigroups with complete intersection cones is 
	the smallest class of affine semigroups closed under the operation of
	s-gluing containing all free affine semigroups. 
\end{theorem} 
  
\begin{remark}  {\rm The notion of p-gluing 
of affine semigroups was introduced in \cite{B-M-T} in order
to characterize the affine toric varieties  
that are set-theoretic complete intersections in positive characteristic $p$. 
It follows from \cite{B-M-T} and  
Theorem 4.3 that any toric variety that is set-theoretic complete intersection 
on binomials in any characteristic 
has a complete intersection cone. \\ 
M. Morales and A. Thoma \cite{M-T} have characterized the subsemigroups $\bn A$ 
with no invertible elements that are complete intersection, $A$ is a finite 
set of elements of $\bz ^n \oplus T$ with $T$ a torsion group. Semigroups of this kind
correspond to lattice ideals.
To any such semigroup, one can associate the 
cone generated by the non torsion parts of the elements of $A$. It follows from the
results in \cite{M-T} and Theorem 4.3 that complete intersection lattice  
ideals have also complete intersection cones.} 
\end{remark}  
{\em Acknowledgment.} The second author would like to thank the  
Department of Mathematics and Statistics  
of the University of Cyprus for the hospitality during the preparation of this work.

\par {\sc Department of Mathematics and Statistics, University of Cyprus,\\  
1678 Nicosia (CYPRUS)}\\ 

{\sc Department of Mathematics, 
University of Ioannina, Ioannina 45110 (GREECE)}  

\begin{thebibliography}{00}  
 
\bibitem{B-M-T} 
	M. Barile, M. Morales and A. Thoma, {\em Set-theoretic 
	complete intersections on binomials}, Proc. Amer. Math. Soc.  {\bf 130} (2002) 1893-1903.  
\bibitem{Ewald}  
	G. Ewald, Combinatorial Convexity and Algebraic Geometry. Graduate Texts  
	in Mathematics {\bf 168},   
	Springer Verlag New York, Berlin, Heidelberg 1996.   
\bibitem{F-M-S}  
	K. Fischer, W. Morris, J. Shapiro, {\em Affine semigroup rings  
	that are complete intersections}, Proc. Amer. Math. Soc.  {\bf 125} (1997) 3137-3145.  
\bibitem{K}  
	M. Katzman, {\em Bipartite graphs whose edge algebras are complete intersection},  
	J. Algebra {\bf 220} (1999) 519-530.  
\bibitem{M-T}  
	M. Morales and A. Thoma, {\em Complete intersection lattice ideals}, preprint 2004.  
\bibitem{R-Z} 
	J. Richter-Gebert and G. M. Ziegler, {\em Oriented matroids}, Handbook on Discrete and Computational Geometry, 
	CRC Press, Boca Raton, FL, (1997) 111-132. 
\bibitem{R-G}   
	J.C. Rosales and P. A. Garcia-Sanchez, {\em On complete intersection  
	affine semigroups},  Commun. Algebra,  {\bf 23}(14) (1995) 5395-5412.  
\bibitem{St}   
	B. Sturmfels, Gr{\"o}bner Bases and Convex Polytopes.  
	University Lecture  Series, No. 8 American Mathematical Society Providence, R.I. 1995.  
 
  
\end{thebibliography}
\end{document}